# TESTING THE SUITABILITY OF POLYNOMIAL MODELS IN ERRORS-IN-VARIABLES PROBLEMS

By Peter Hall and Yanyuan Ma

*Australian National University, and University of Neuchatel and Texas A&M University*

A low-degree polynomial model for a response curve is used commonly in practice. It generally incorporates a linear or quadratic function of the covariate. In this paper we suggest methods for testing the goodness of fit of a general polynomial model when there are errors in the covariates. There, the true covariates are not directly observed, and conventional bootstrap methods for testing are not applicable. We develop a new approach, in which deconvolution methods are used to estimate the distribution of the covariates under the null hypothesis, and a "wild" or moment-matching bootstrap argument is employed to estimate the distribution of the experimental errors (distinct from the distribution of the errors in covariates). Most of our attention is directed at the case where the distribution of the errors in covariates is known, although we also discuss methods for estimation and testing when the covariate error distribution is estimated. No assumptions are made about the distribution of experimental error, and, in particular, we depart substantially from conventional parametric models for errors-in-variables problems.

**1. Introduction.** Suppose we observe independent pairs $(W_1, Y_1), \ldots, (W_n, Y_n)$ distributed as $(W, Y)$, where

(1) $\qquad Y = g(X) + \varepsilon, \qquad W = X + U, \qquad E(\varepsilon \,|\, U, X) = 0,$

and $U$ is independent of $X$ and has zero mean. The particular model of interest is that where $g$ is a polynomial,

(2) $$g(x) = \sum_{j=0}^{p} \beta_j^0 x^j.$$

Here, for $0 \le j \le p$, $\beta_j^0$ denotes the true value of a parameter $\beta_j$.

Our main purpose in this paper is to suggest ways of assessing the goodness of fit of the polynomial model. We shall treat the goodness-of-fit problem as one of testing the null hypothesis that $g(x)$ can be expressed as in (2), for values of $x$ that lie in the support of the distribution of $X$. Note that in the absence of measurement error, that is, if we could observe $X$, then this problem could be solved readily by using, for example, the test statistic and its properties developed by Fan, Zhang and Zhang [14]. However, the presence of the measurement error complicates the problem, and we are not aware of an existing method in that case.

A related problem has been treated by Cheng and Kukush [4]. There, an ingenious, asymptotic squared-difference goodness-of-fit test is suggested, based on a statistic which, under the null hypothesis, has a limiting chi-squared distribution with one degree of freedom. However, it is readily seen that, while the Cheng and Kukush [4] test has good power properties against some alternatives, it has zero power against many others. Intuitively, this is because the test addresses only one mode of a potentially infinite number of modes of departure from the null hypothesis of a polynomial fit. That single mode, or single component in the infinite class of components that are all orthogonal to the class of all polynomials of degree $p$, is responsible for the single degree of freedom in the test of Cheng and Kukush [4].

By way of contrast, the test proposed in the present paper addresses simultaneously the infinity of components that can define departure from the class of polynomials of degree $p$. In this setting the limiting distribution of any test statistic will be relatively complex, and an asymptotic test will not be feasible. We suggest instead a bootstrap method for calibrating the test and producing critical points.

However, bootstrap methods in this problem are necessarily quite non-standard. Indeed, the bootstrap is seldom used in the context of errors in variables, since neither the explanatory variables $X$ nor the errors $\varepsilon$ can be directly accessed. At best, only their distributions can be estimated, and so the bootstrap cannot proceed by resampling either observed or imputed data, such as residuals.

Quite different methods are required for estimating the distributions of $X$ and $\varepsilon$, as a prelude to applying the bootstrap. From some points of view, estimating the distribution of $X$ is the simpler of the two tasks; that problem is one of conventional deconvolution, in which, given the distribution of $U$, we wish to estimate the distribution of $X$ from data on $W = X + U$. However, it can be shown theoretically that, unless the distribution of $U$ is especially unsmooth, the distribution function of $X$ cannot be estimated root-$n$ consistently. For example, if the characteristic function of the distribution of $U$ decays like $|t|^{-\alpha}$ as $|t| \to \infty$, where $\alpha > 0$, then it can be proved that a necessary condition for root-$n$ consistency to be achievable is that $\alpha < \frac{1}{2}$. This constraint denies even a single derivative to the density of $U$.

Therefore, the distribution of $X$ seldom can be estimated root-$n$ consistently, and so bootstrapping that variable presents challenges. Estimating the distribution of $\varepsilon$ is an even more awkward problem. However, a careful examination of theoretical issues shows that the limiting distribution of our test statistic depends on properties of $\varepsilon$ only to the extent of $\operatorname{var} \varepsilon$, and so the moment-matching, or wild, bootstrap is feasible for estimating the distribution of experimental error.

In summary, key contributions of this paper are constructing a statistic for testing model adequacy in the context of polynomial models with measurement error; proposing a nonstandard bootstrap method for assessing the distribution of the test statistic under the null hypothesis; and showing how to use repeated measurements when the measurement-error distribution is unknown. Innovations include the novel form of the test statistic, the unconventional way in which it is computed, using both deconvolution and wild-bootstrap techniques, and our theoretical derivation of properties of the test.

Various forms of (1) have been studied in the literature. Most of the work focuses on a parametric model framework, where a parametric form of the distribution of $\varepsilon$ given $U$ and $X$ is adopted, typically being normal. When $g(X)$ is linear, extensive research can be found in Fuller [17], and the efficient estimator was given by Bickel and Ritov [1]. The same efficient estimator was also discovered in a broader generalized linear model framework by Stefanski and Carroll [27]. The extension of $g(X)$ to a general polynomial was first studied in Chan and Mak [3], where a root-$n$ consistent estimator was constructed. Their work was later further extended by Cheng and Schneeweiss [5] and Cheng, Schneeweiss and Thamerus [6]. A comparison of several methods is given by Kukush, Schneeweiss and Wolf [21].

A review and study of a class of estimators can be found in Taupin [29]. Consistent and efficient estimators for a general function $g$ were recently constructed by Tsiatis and Ma [30]. Estimators proposed in Bickel and Ritov [1] and Cheng and Schneeweiss [5] also apply when a distributional model for $(\varepsilon|U, X)$ is not assumed, hence their model is in fact semiparametric. A further extension from this semiparametric model framework is to consider a partially linear model through replacing $g(X)$ by $X\beta + \theta(Z)$, where $\beta$ is an unknown parameter and $\theta(Z)$ is an arbitrary unknown function of some observable covariates $Z$. Estimators in this setting were proposed by Liang, Härdle and Carroll [23]. When no functional form is assumed for $g(X)$, the model becomes nonparametric. Estimators and their properties are studied in Fan and Truong [13] and Efromovich [10, 11]. Recent work on the moment-matching bootstrap includes that of Fan and Li [15], Flachaire [16], Domínguez and Lobato [9], Prášková [26], Kauermann and Opsomer [20], Li, Hsiao and Zinn [22] and González Manteiga, Martínez Miranda and Pérez González [18].

The distribution of $U$ in (1) might be known, or it might be estimated directly or from replicated data on $W$. To focus on the main problem, we assume first that the distribution of $U$ is known, and treat subsequently, in Section 5, the case where it is unknown.

## 2. Methodology.

2.1. *Methodology for estimating $\beta_0, \ldots, \beta_p$.* Because various methods for estimating $\beta$ already exist, we only briefly outline the estimator that is used in this paper. Let $x_0, x_1, \ldots$ be real numbers with $x_0 \neq 0$, and define recursively functions $P_j$ of $j+1$ variables by $P_0(x_0) = x_0^{-1}$, and

$$x_0 P_j(x_0, \ldots, x_j) = -\sum_{k=0}^{j-1} \binom{j}{k} x_{j-k} P_k(x_0, \ldots, x_k), \qquad j \geq 1.$$

Given a random variable $Z$, define $\mu_j(Z) = E(Z^j)$. From the data $(W_k, Y_k)$, construct estimators of $a_j = \mu_j(W)$ and $A_j = E(YW^j)$ using $\hat{a}_j = n^{-1} \sum_k W_k^j$ and $\hat{A}_j = n^{-1} \sum_k Y_k W_k^j$, respectively. Define $b_j = \mu_j(X)$ and $B_j = E(YX^j)$, and put $\nu_j = \mu_j(U)$, a known quantity. It can be shown that

$$b_j = \sum_{k=0}^{j} \binom{j}{k} a_{j-k} P_k(\nu_0, \ldots, \nu_k),$$

$$B_j = \sum_{k=0}^{j} \binom{j}{k} A_{j-k} P_k(\nu_0, \ldots, \nu_k).$$

Hence, under moment assumptions, root-$n$ consistent estimators of $b_j$ and $B_j$ are

(3)
$$\hat{b}_j = \sum_{k=0}^{j} \binom{j}{k} \hat{a}_{j-k} P_k(\nu_0, \ldots, \nu_k),$$

$$\hat{B}_j = \sum_{k=0}^{j} \binom{j}{k} \hat{A}_{j-k} P_k(\nu_0, \ldots, \nu_k).$$

An estimator of the true values $\beta^0 = (\beta_0, \ldots, \beta_p)^{\mathrm{T}}$ is given by

(4) $$\hat{\beta} = \widehat{M}^{-1} \hat{B},$$

where $\hat{B} = (\hat{B}_0, \ldots, \hat{B}_p)^{\mathrm{T}}$, $\widehat{M} = (\hat{m}_{jk})$ is a $(p+1) \times (p+1)$ matrix, and $\hat{m}_{jk} = \hat{b}_{j+k}$ for $0 \leq j, k \leq p$. It can be proved that $\hat{\beta}$ is root-$n$ consistent for $\beta^0$ and is asymptotically normally distributed, provided (1) and (2) hold, $E(W^{4p}) + E(\varepsilon^2) < \infty$ and the distribution of $X$ has a nondegenerate continuous component.

2.2. *Hypotheses and test statistic.* Consider the problem of testing the null hypothesis $H_0 = H_0(p)$, that $g$ in the model (1) is given by (2) for an appropriate choice of $\beta_0, \ldots, \beta_p$, against the complementary alternative $H_1(p)$. Since we have access to information about $g(x)$ only when $x$ is in the support of the density $f_X$ of $X$, then $H_1(p)$ should have the form: $g$ is not equal almost everywhere on $\mathrm{supp} f_X$ to a polynomial of degree $p$. Equivalently, $H_1(p)$ is characterized by the class of functions $g$ that are bounded on compact intervals and satisfy

$$\inf_{\beta_0,\ldots,\beta_p} \int_{-t_0}^{t_0} |\psi(t) - \phi(t\,|\,\beta)|^2 \, dt > 0$$

for each $t_0 > 0$, where

$$\psi(t) = E\{g(X)e^{itX}\}, \qquad \phi(t\,|\,\beta) = E\left(e^{itX} \sum_{j=0}^{p} \beta_j X^j\right)$$

and $i = \sqrt{-1}$. In defining $\psi$ and $\phi$ we assume that $E\{|g(X)| + |X|^p\} < \infty$.

These considerations suggest that we base our test on the statistic $T(\tilde{\beta})$, where

(5) $$T(\beta) = \int |\hat{\psi}(t) - \hat{\phi}(t\,|\,\beta)|^2 w_1(t) \, dt;$$

$\hat{\psi}(t)$ and $\hat{\phi}(t\,|\,\beta)$ are root-$n$ consistent estimators of $\psi(t)$ and $\phi(t\,|\,\beta)$, respectively; $\tilde{\beta}$ denotes either $\hat{\beta}$, defined at (4), or an alternative estimator, such as $\mathrm{argmin}\, T(\beta)$; and $w_1 > 0$ is a known weight function.

For computational simplicity, we shall take $\tilde{\beta} = \hat{\beta}$. Unbiased estimators of $\psi(t)$ and $\phi(t\,|\,\beta)$ are given by

(6)
$$\hat{\psi}(t) = \frac{1}{n f_U^{\mathrm{Ft}}(t)} \sum_{j=1}^{n} Y_j e^{itW_j},$$

$$\hat{\phi}(t\,|\,\beta) = \left\{\sum_{k=0}^{p} \beta_k (i^{-1} D_t)^k\right\} \left\{\frac{\sum_j e^{itW_j}}{n f_U^{\mathrm{Ft}}(t)}\right\},$$

where $D_t = \partial/\partial t$ is the differentiation operator, and $f_U^{\mathrm{Ft}}$ is the characteristic function of $U$, or equivalently, the Fourier transform of $f_U$. For these choices of $\hat{\beta}$, $\hat{\psi}$ and $\hat{\phi}$, our test amounts to rejecting $H_0$ if the statistic $S \equiv T(\hat{\beta})$ is too large. We shall use bootstrap methods to determine a critical point for the test. As a prelude to that step, we require estimators of the distributions of $X$ and $\varepsilon$.

In order to remove the function $f_U^{\mathrm{Ft}}$ from denominators in (6), it is convenient to take $w_1 = (f_U^{\mathrm{Ft}})^2 w$, where $w$ is another weight function. This produces the statistic

$$T(\beta) = \int |\hat{\psi}(t) - \hat{\phi}(t\,|\,\beta)|^2 f_U^{\mathrm{Ft}}(t)^2 w(t) \, dt.$$

2.3. *Estimator of distribution of $X$.* The distribution of $X$ is accessible using conventional deconvolution methods, as follows. Given data $W_1, \ldots, W_n$ on $W = X + U$, a kernel estimator of the density $f_X$ of $X$ is given by

$$\tilde{f}_X(x) = \tilde{f}_X(x \mid h) = \frac{1}{nh} \sum_{j=1}^{n} L\left(\frac{x - W_j}{h}\right),$$

where

$$L(u) = \frac{1}{2\pi} \int e^{-itu} \frac{K^{\mathrm{Ft}}(t)}{f_U^{\mathrm{Ft}}(t/h)} \, dt,$$

$K$ is a kernel function (in particular, a function which integrates to 1), $K^{\mathrm{Ft}}$ denotes the Fourier transform of $K$ and $h > 0$ is a smoothing parameter. See, for example, Carroll and Hall [2], Stefanski and Carroll [28] and Fan [12].

Integrating $\tilde{f}_X$, we obtain an estimator $\widetilde{F}_X$ of the distribution function $F_X$ of $X$,

$$\widetilde{F}_X(x) = \widetilde{F}_X(x \mid h) = \int_{-\infty}^{x} \tilde{f}_X(u) \, du = \frac{1}{n} \sum_{j=1}^{n} L_1(x - W_j),$$

where $L_1(hu) = \int_{v \leq u} L(v) \, dv$, or equivalently,

(7) $$L_1(u) = \frac{1}{2} + \frac{1}{2\pi} \int_{-\infty}^{\infty} \frac{\sin tu}{t} \frac{K^{\mathrm{Ft}}(ht)}{f_U^{\mathrm{Ft}}(t)} \, dt,$$

provided that $K^{\mathrm{Ft}}(ht)/f_U^{\mathrm{Ft}}(t)$ is real valued. In Section 3 we discuss choice of $K$ and $h$.

Next we convert $\widetilde{F}_X$ to a distribution function $\hat{F}_X$ by defining first $\bar{F}_X(x) = \max_{u \leq x} \widetilde{F}_X(u)$ and then

(8) $$\hat{F}_X(x) = \frac{\bar{F}_X(x) - \bar{F}_X(c_1)}{\bar{F}_X(c_2) - \bar{F}_X(c_1)}$$

if $c_1 \leq x < c_2$, $\hat{F}_X(x) = 0$ if $x < c_1$, and $\hat{F}_X(x) = 1$ if $x \geq c_2$, where $c_1 < c_2$ are constants.

2.4. *Estimator of distribution of $\varepsilon$.* Conventional deconvolution methods can be used to estimate the distribution of $\varepsilon$ when $p = 1$, although they are awkward to implement; and they fail for $p > 1$. Fortunately, satisfactory accuracy can be obtained using a simpler, moment-matching or "wild" bootstrap approach. To this end, let $\omega_r = E(\varepsilon^r)$, for integers $r \geq 1$, and note that $\omega_1 = 0$; let $\hat{\omega}_r$, for $r \geq 2$, denote respective estimators of $\omega_r$;

let $G(\cdot\,|\,\kappa_2,\ldots,\kappa_q)$ be a known distribution with zero mean and moments $\int x^r\,dG(x) = \kappa_r$, for $2 \leq r \leq q$, where $q \geq 2$; and put

$$\hat{F}_\varepsilon = G(\cdot\,|\,\widehat{\omega}_2,\ldots,\widehat{\omega}_q). \tag{9}$$

This estimator is generally not consistent for $F_\varepsilon$, but it is adequate for our purpose. Examples of the distribution $G$ will be given in Section 3.

In Section 4 we shall show that the asymptotically correct level for the test is achieved by taking $q \geq 2$. Although $q = 2$ is sufficient, accuracy can be improved by using $q = 3$, or, in the case of near symmetry, fitting a distribution with first and third moments equal to zero and second and fourth moments equal to $\widehat{\omega}_2$ and $\widehat{\omega}_4$; see Section 3. The moment-matching bootstrap could also be employed to estimate the distribution of $X$, but there we require $q \geq 4p$.

Next we define estimators $\widehat{\omega}_r$. Observe that

$$\begin{aligned}
\omega_r = E(\varepsilon^r) &= E(Y^r) - \sum_{s=0}^{r-1} \binom{r}{s} E(\varepsilon^s) E\left(\sum_{j=0}^{p} \beta_j X^j\right)^{r-s} \\
&= E(Y^r) - \sum_{s=0}^{r-1} \sum_{t_0+\cdots+t_p=r-s} \frac{r!}{s!\,t_0!\cdots t_p!} \omega_s \beta_0^{t_0}\cdots \beta_p^{t_p} E(X^{t_1+2t_2+\cdots+pt_p}),
\end{aligned} \tag{10}$$

where the second summation is over integers $t_0,\ldots,t_p \geq 0$ such that $t_0+\cdots+t_p = r-s$, and, since $E(\varepsilon) = 0$, we may exclude from the first summation in (10) the term corresponding to $s = 1$. Results (3) and (4) give us root-$n$ consistent estimators $\hat{b}_j$ and $\hat{\beta}_j$ of $b_j = E(X^j)$ and $\beta_j$, respectively, and we can readily compute $\bar{Y}_r = n^{-1}\sum_j Y_j^r$, an unbiased estimator of $E(Y^r)$. Therefore, having constructed estimators $\widehat{\omega}_1 = 0$, $\widehat{\omega}_2,\ldots,\widehat{\omega}_{r-1}$, we define $\widehat{\omega}_r$ recursively by

$$\widehat{\omega}_r = \bar{Y}_r - \sum_{s=0}^{r-1} \sum_{t_0+\cdots+t_p=r-s} \frac{r!}{s!\,t_0!\cdots t_p!} \widehat{\omega}_s \hat{\beta}_0^{t_0}\cdots \hat{\beta}_p^{t_p} \hat{b}_{t_1+2t_2+\cdots+pt_p}.$$

2.5. *Implementing the bootstrap test.* Our bootstrap method has six steps, as follows. (a) Compute the estimators $\hat{\beta}_0,\ldots,\hat{\beta}_p$ suggested in Section 2.1, and the distribution estimators $\hat{F}_X$ and $\hat{F}_\varepsilon$ suggested at (8) and (9). Calculate the test statistic $S = T(\hat{\beta})$ from (5). (b) Draw data $X_1^*,\ldots,X_n^*$ from $\hat{F}_X$, $\varepsilon_1^*,\ldots,\varepsilon_n^*$ from $\hat{F}_\varepsilon$ and $U_1^*,\ldots,U_n^*$ from the distribution of $U$, and put $\hat{g}(x) = \sum_{0 \leq j \leq p} \hat{\beta}_j x^j$, $Y_j^* = \hat{g}(X_j^*) + \varepsilon_j^*$ and $W_j^* = X_j^* + U_j^*$, for $1 \leq j \leq n$. (c) Using the data pairs $(W_j^*, Y_j^*)$ in place of $(W_j, Y_j)$, compute the estimator $\hat{\beta}^* = (\hat{\beta}_0^*,\ldots,\hat{\beta}_p^*)^\mathrm{T}$ of $\beta = (\beta_0,\ldots,\beta_p)^\mathrm{T}$. (d) Compute the analogue $T^*(\beta)$ of $T(\beta)$ defined at (5), using $(W_j^*, Y_j^*)$ instead of $(W_j, Y_j)$, and form the statistic $S^* = T^*(\hat{\beta}^*)$, the bootstrap analogue of $S = T(\hat{\beta})$. (e) Using repeated

Monte Carlo simulation, approximate the distribution of $S^*$ conditional on the data $\mathcal{D} = \{(W_1, Y_1), \ldots, (W_n, Y_n)\}$, and in particular, approximate the critical point $\hat{s}_\alpha$ such that $P(S^* > \hat{s}_\alpha | \mathcal{D}) = 1 - \alpha$. (f) Reject $H_0(p)$ in favor of $H_1(p)$ at the nominal level $\alpha$ if $S > \hat{s}_\alpha$.

## 3. Computational issues and numerical results.

3.1. *Choice of $K$ and $h$.* Generally, $K$ is selected so that $K^{\text{Ft}}$ vanishes outside a compact interval. A popular choice is

$$(11) \quad K(x) = 48(\cos x)(1 - 15x^{-2})(\pi x^4)^{-1} - 144(\sin x)(2 - 5x^{-2})(\pi x^5)^{-1},$$

for which $K^{\text{Ft}}(t) = (1 - t^2)^3$ if $|t| \leq 1$ and $K^{\text{Ft}}(t) = 0$ otherwise. See, for example, Delaigle and Gijbels [7, 8].

In such cases, $L_1(u)$ is well defined by (7) and finite for each $u$, provided

$$(12) \qquad f_U^{\text{Ft}} \text{ is real-valued and does not vanish on the real line.}$$

This is a common assumption in deconvolution problems, and while it can be circumvented, we shall use models for which it holds.

As a prelude to bandwidth choice, we note that if (12) holds and $f_U^{\text{Ft}}(t) \sim Ct^{-\alpha}$ as $t \to \infty$, with $\alpha > \frac{1}{2}$, and if $K$ is as at (11), then

$$\int_{-\infty}^\infty E\{\widetilde{F}_X(x \,|\, h) - F_X(x)\}^2 \, dx = C_1 n^{-1} h^{1-2\alpha} + C_2 h^4 + o(n^{-1} h^{1-2\alpha} + h^4),$$

where $C_1 = C^2 \kappa(\alpha)/\pi$, $C_2 = 9 \int (f_X')^2 \, dx$ and $\kappa(\alpha) = \int_{t>0} t^{2\alpha-2} K^{\text{Ft}}(t)^2 \, dt$.

Delaigle and Gijbels [8] suggested methods for estimating $J_X = \int (f_X')^2 \, dx$, and hence, for approximating $C_2$. The simplest of their techniques is a "normal reference" approach, analogous to bandwidth choice in density estimation by comparison with the normal distribution. Specifically, $f_X$ is taken to be a normal $N(0, \sigma_X^2)$ density, where $\sigma_X^2 = \text{var } X$ and is estimated by $\hat{\sigma}_X^2$, equal to the empirical variance of the data $W_1, \ldots, W_n$, minus the known variance of $U$. If $X$ were normally distributed, then $J_X$ would equal $(4\pi^{1/2} \sigma_X^3)^{-1}$. Therefore, we take $\hat{J}_X = (4\pi^{1/2} \hat{\sigma}_X^3)^{-1}$, and so our estimator of $C_2$ is $\hat{C}_2 = 9/(4\pi^{1/2} \sigma_X^3)$. Finally, we compute $\kappa(\alpha)$, and then $C_1$, using the known value of $\alpha$, and choose $h$ to minimize $C_1 n^{-1} h^{1-2\alpha} + \hat{C}_2 h^4$.

3.2. *Choice of the distribution $G = G(\cdot \,|\, \omega_2, \ldots, \omega_q)$.* The simplest case is that where $q = 2$, in which instance one would generally take $G$ to be the normal $N(0, \omega_2)$ distribution. Two examples of distributions $G$ that are suitable when $\omega_1 = \omega_3 = 0$ and $q = 4$ are the three-point distribution defined by

$$(13) \qquad P(Z = 0) = 1 - \pi, \qquad P(Z = \pm \pi^{-1/2}) = \tfrac{1}{2}\pi,$$

where $0 < \pi < 1$, and the Student's $t$ distribution. The three-point distribution can be used to capture any pair $(\omega_2, \omega_4)$, regardless of the sign of kurtosis. The Student's $t$ distribution can capture only $(\omega_2, \omega_4)$ for which kurtosis is positive; however, the positive sign is the more common in practice.

The distribution at (13) has $E(Z) = 0$, $E(Z^2) = 1$ and $E(Z^4) = \pi^{-1}$. Therefore, if $\pi = \omega_2^2/\omega_4$, then the distribution of $\omega_2^{1/2} Z$ is symmetric with variance and fourth moment equal to $\omega_2$ and $\omega_4$, respectively. To implement the moment-matching method in this setting, one replaces $\omega_2$ and $\omega_4$ by their respective estimates, $\widehat{\omega}_2$ and $\widehat{\omega}_4$, discussed in Section 2.4; takes the estimator of $\pi$ to be $\widehat{\omega}_2^2/\widehat{\omega}_4$; estimates the distribution of $Z$, at (13), by replacing $\pi$ by this estimator; and, if $Z$ has this estimated distribution, takes the distribution $G$ (our surrogate for the distribution of $\varepsilon$) to be that of $\widehat{\omega}^{1/2} Z$. In some instances it is not possible to reliably estimate high-order moments, and there only low-order moments are fitted. For example, in the Alaskan Earthquake example in Section 3.4 we fit the normal $N(0, \omega_2)$ distribution.

The convergence rate of the wild bootstrap typically improves when higher moments are fitted in additional to the first two moments, as discussed by Liu [24], Mammen [25] and Härdle and Mammen [19]. The model that they consider has the form $Y_i = g(X_i) + \varepsilon_i$, where $\varepsilon_i$ is not necessarily identically distributed. When $g$ is linear, Liu [24] shows analytically that the second-order properties of the wild bootstrap are obtained when the third moment is fitted, due to a correction of a skewness term in the Edgeworth expansion of a sampling distribution. Following these results, in the following simulation studies we implement the three-point distribution which fits the first four moments.

3.3. *Simulation results.* We conduct two simulation studies. In the first, we generate data from the linear errors-in-variables model $Y_i = \beta_1 X_i + \beta_0 + \varepsilon_i$, with $W_i = X_i + U_i$, $i = 1, \ldots, n$. The latent variables $X_i$ are generated from a uniform distribution in $[-3, 4]$, and the experimental errors $\varepsilon_i$ come from a normal distribution with mean 0 and variance 1. We generate the measurement errors $U_i$ from two different distributions: a normal $N(0, 1)$ distribution and a Laplace distribution with variance 0.5. In each parameter setting we simulate 2000 datasets, for various sample sizes $n$, and use bootstrap resample size 100. In the first simulation study, datasets are generated under $H_0$.

The purpose here is to assess the level accuracy of the test. Results are given in Table 1. The two different measurement error distributions in the upper and lower halves of the table correspond to two approaches to choosing $h$. In the case of normal error, the optimal $h$ is infinity, while for Laplace error, a finite value of $h$ is obtained using the strategy described in Section 3.1. As can be seen from Table 1, even in this simple model and even for very small sample sizes $n = 50$ to $n = 100$, the rejection levels under the null hypothesis are close to the desired levels for both error types. We repeated the simulations with bootstrap resample size 200 and obtained very similar results (not reported here).

TABLE 1
*Simulation 1: Level accuracy. In the first block (upper half table), the measurement error $U_i$ is normal $N(0,1)$ and the bandwidth is $h = 5.0$, which is practically infinite with respect to the support of the distribution of $X$. In the second block (lower half table), $U_i$ has a Laplace distribution with zero mean and variance $0.5$, and the bandwidth $h$ is calculated as suggested in Section 3.1. Each entry in the table is based on 2000 simulated datasets. The bootstrap resample size is 100. The very top row of the table gives the nominal levels*

| $n$ | 5% | 6% | 7% | 8% | 9% | 10% |
|-----|-----|-----|-----|-----|-----|-----|
| | | | Normal measurement error | | | |
| 50 | 4.55% | 5.85% | 6.38% | 7.63% | 9.45% | 9.93% |
| 60 | 4.38% | 6.05% | 6.68% | 8.00% | 10.00% | 10.50% |
| 70 | 4.53% | 6.05% | 6.55% | 7.63% | 9.10% | 9.63% |
| 80 | 5.03% | 6.15% | 6.70% | 7.75% | 9.20% | 9.73% |
| 90 | 4.33% | 5.60% | 6.13% | 7.28% | 9.20% | 9.70% |
| 100 | 5.33% | 6.60% | 6.95% | 7.93% | 9.65% | 10.13% |
| | | | Laplace measurement error | | | |
| 50 | 4.75% | 6.40% | 6.83% | 7.75% | 9.30% | 10.05% |
| 60 | 5.05% | 6.65% | 7.23% | 8.45% | 10.35% | 10.88% |
| 70 | 4.78% | 6.55% | 7.13% | 8.23% | 9.55% | 10.08% |
| 80 | 4.45% | 6.30% | 6.83% | 7.95% | 9.65% | 10.20% |
| 90 | 5.25% | 6.50% | 7.18% | 8.50% | 9.75% | 10.40% |
| 100 | 5.05% | 6.70% | 7.15% | 8.15% | 9.90% | 10.65% |

The second simulation study addresses power. Here we generate datasets from the model $Y_i = \beta_1 X_i + \beta_0 + c\cos(X_i) + \varepsilon_i$, with $W_i = X_i + U_i$ and $c = 1.5$ a constant. (The first study used the same model but with $c = 0$. The distributions of $\varepsilon_i$, $U_i$ and $X_i$ are as in the first study.) We again take bootstrap resample size to be 100, but this time we calculate the power of the test at different levels. The results for different sample sizes are given in Table 2. We can see that as sample size increases, so too does the power. In this experiment, sample size $n = 80$ can already achieve 90% power at level 5%. For $n = 100$, power increases to 95%, which is usually sufficient in practice. In simulations with bootstrap resample size 200 we obtained very similar results. Hence, bootstrap resample size 100 seems adequate.

3.4. *Alaskan Earthquake example.* We implement our method on the Alaskan Earthquake data studied by Fuller ([17], Chapters 1 and 4). In this

TESTING ERRORS IN VARIABLES MODELS 11TABLE 2
*Simulation 2: Power. Parameter settings are as for Table 1*

| $n$ | 5% | 6% | 7% | 8% | 9% | 10% |
|-----|-----|-----|-----|-----|-----|-----|
| | | Normal measurement error | | | | |
| 50 | 71.98% | 75.90% | 76.65% | 78.05% | 79.80% | 80.40% |
| 60 | 78.78% | 81.60% | 82.65% | 84.20% | 85.75% | 86.30% |
| 70 | 84.80% | 87.25% | 87.93% | 89.08% | 90.15% | 90.68% |
| 80 | 90.63% | 92.30% | 92.80% | 93.68% | 94.60% | 95.03% |
| 90 | 93.15% | 94.60% | 94.90% | 95.48% | 96.30% | 96.48% |
| 100 | 95.93% | 96.85% | 96.95% | 97.15% | 97.55% | 97.73% |
| | | Laplace measurement error | | | | |
| 50 | 76.93% | 79.80% | 80.70% | 82.40% | 84.65% | 85.28% |
| 60 | 85.55% | 87.85% | 88.50% | 89.78% | 91.00% | 91.35% |
| 70 | 91.08% | 92.55% | 93.03% | 93.73% | 94.45% | 94.70% |
| 80 | 93.88% | 95.05% | 95.33% | 95.80% | 96.50% | 96.58% |
| 90 | 95.80% | 96.70% | 97.03% | 97.58% | 98.05% | 98.10% |
| 100 | 97.10% | 97.95% | 98.25% | 98.68% | 99.00% | 99.10% |

dataset, the logarithm of the seismogram amplitude of 20 second surface waves ($Y$) and the logarithm of the seismogram amplitude of longitudinal body waves ($W$) of 62 earthquakes are recorded. The main interest is in analyzing how the surface wave is related to the longitudinal body wave. Of course, both variables are measured with error. Fuller [17] used a linear errors-in-variable model, $Y = \beta_0 + \beta_1 X + \varepsilon$ and $W = X + U$. Assuming a normal $N(0, \sigma_U^2)$ for $U$, and using extra available information, Fuller [17] estimated the measurement error variance to be $\sigma_U^2 = 0.035$, with standard error 0.0086.

We implement a wild bootstrap procedure that matches the first three estimated moments of the experimental error distribution. (The fourth moment estimate here is too highly variable to be reliable, and, in fact, its point estimate is negative.) Taking $\sigma_U^2 = 0.035$, and employing bootstrap resample sizes 100, 200, 300, 400 or 500, we find that the resulting $p$-values are all between 54% and 57%.

Considering that the value of $\sigma_U^2$ is estimated, we also consider two extreme cases, where $\sigma_U^2$ equals $0.035 \pm 3.5 \times 0.0086 = 0.0049$ and 0.065, respectively, and apply the testing procedure in these cases. The results associated with these values of $\sigma_U^2$ and with different bootstrap resample sizes are reported in Table 3.

For none of these parameter settings is the reported $p$-value small enough to cast reasonable doubt on the adequacy of the linear model for this dataset. Although the different values of $\sigma_U^2$ cause significant changes to estimators of $\beta_0$ and $\beta_1$, the evidence in favor of rejecting the null hypothesis is virtually nonexistent, in all three cases. This observation reflects the substantial

TABLE 3
*Alaskan Earthquake example: Testing for linearity between the logarithm of the seismogram amplitude of 20 second surface waves (Y), and the logarithm of the seismogram amplitude of longitudinal body waves (X). Three different values for the measurement error variance $\sigma_U^2$ are considered. The p-value of the test, $p(B)$, is reported as a function of the number of bootstrap resamples, B*

| $\sigma_U^2$ | $\hat{\beta}_0$ | $\hat{\beta}_1$ | $p(100)$ | $p(200)$ | $p(300)$ | $p(400)$ | $p(500)$ |
|---|---|---|---|---|---|---|---|
| 0.0049 | $-1.65$ | 1.29 | 50.0% | 44.0% | 43.3% | 45.0% | 45.4% |
| 0.035 | $-2.81$ | 1.51 | 55.0% | 55.5% | 54.0% | 54.5% | 56.8% |
| 0.065 | $-4.47$ | 1.83 | 70.0% | 72.5% | 72.3% | 73.5% | 75.8% |

variability in the dataset, noted from a different viewpoint two paragraphs above.

## 4. Theoretical properties.
We shall assume that

(14) $E(W^{4p}) < \infty$, $0 < E(\varepsilon^2) < \infty$, the distribution of $X$ has a nondegenerate continuous component and $f_U^{\text{Ft}}$ vanishes only at isolated points;

(15) $w(t) > 0$ for each $t$, $w(t)$ converges to zero faster than any polynomial as $|t| \to \infty$ and $\max_{1 \le k \le p} \int |D_t^k f_U^{\text{Ft}}(t)^{-1}|^2 f_U^{\text{Ft}}(t)^2 w(t)\, dt < \infty$.

In the context of conventional models for the distribution of $U$, $|D_t^k f_U^{\text{Ft}}(t)|/f_U^{\text{Ft}}(t)$ is dominated by a polynomial in $t$, and in such cases the last part of (15) follows from the rest of that assumption.

When implementing the bootstrap we shall assume, in addition to (14) and (15), that

(16) the support of the distribution of $X$ is contained within the finite interval $(c_1, c_2)$,

where $c_1$ and $c_2$ are fixed and are used in the definition of $\hat{F}_X$ at (8).

Of course, the distribution $G = G(\cdot \mid \kappa_2, \ldots, \kappa_q)$ that we employ to estimate $F_\varepsilon$, at (9), has by definition finite variance if $\kappa_2 < \infty$, so we do not impose this as a regularity condition. It is not necessary to stipulate whether the distribution $G$ is discrete or continuous.

The main theoretical properties of our estimator are given in the following theorem. There, part (a) describes limit theory under the null hypothesis $H_0(p)$, part (b) asserts consistency of the bootstrap estimator of the distribution of the test statistic under $H_0(p)$, and part (c) shows that the test is able to detect a large class of semiparametric, root-$n$ departures from the null hypothesis. It is straightforward to prove a version of part (b) when $H_0(p)$ fails; that result requires conditions on a class of $g$'s for which $H_1(p)$ holds.

THEOREM 1. *Assume that the data on which the test is based are generated by the model* (1). *Then:* (a) *If* (14) *and* (15) *hold, and if the null hypothesis* $H_0(p)$ *is valid [i.e., if* (2) *holds for some choice of the parameters* $\beta_0^0, \ldots, \beta_p^0$], *then* $nT(\hat{\beta}) \to \xi$ *in distribution, where* $\xi$ *denotes a random variable for which* $P(0 < \xi < \infty) = 1$, *and the distribution of* $\xi$ *depends on that of* $\varepsilon$ *only through* $\mathrm{var}\,\varepsilon$. (b) *If* (14)–(16) *hold, and* $H_0(p)$ *is valid, then the distribution of* $T^*(\hat{\beta}^*)$, *conditional on the data, converges in probability to that of* $\xi$. (c) *If* (14)–(16) *hold, and if the function* $g = g_n$ *in* (1) *is taken to depend on* $n$, *as*

$$(17) \qquad g(x) = \sum_{j=0}^{p} \beta_j^0 x^j + n^{-1/2} c\gamma(x),$$

*where* $\beta_0^0, \ldots, \beta_p^0$ *and* $c > 0$ *are fixed, and the function* $\gamma$ *is bounded, compactly supported and, on a subset of the support of the distribution of* $X$ *that has nonzero Lebesgue measure, does not vanish and does not equal almost everywhere a polynomial of degree* $p$, *then*

$$(18) \qquad \lim_{c \to \infty} \liminf_{n \to \infty} P(S > \hat{s}_\alpha) = 1.$$

The property stated in part (a) of the theorem that the distribution of $\xi$ depends on that of $\varepsilon$ only through $\mathrm{var}\,\varepsilon$, is the key to the fact that the moment-matching bootstrap is adequate for estimating the distribution of $\varepsilon$ when calibrating the test statistic $T(\hat{\beta})$. By way of comparison, the distribution of $\xi$ depends on that of $X$ through more than just the first two or three moments.

The function $g$ at (17) represents a local departure from the null hypothesis $H_0(p)$. Indeed, under the latter hypothesis, $g$ would equal just the first part of the right-hand side of (17). Result (18) asserts that, in the case of a local departure of this form, the test is asymptotically capable of detecting the fact that $H_0(p)$ fails. More particularly, for all sufficiently large $n$, the probability that the test correctly detects the fact that $H_0(p)$ is violated exceeds $1 - \eta$, where $\eta > 0$ can be chosen arbitrarily small by selecting $c$ in (17) sufficiently large.

The assumption that the distribution of $X$ is compactly supported, used in parts (b) and (c) of the theorem, is imposed for convenience and can be relaxed; we do not do so since we wish to keep the proof and the regularity conditions simple.

An outline proof of Theorem 1 will be given in the Appendix. There the distribution of $\xi$ will be given.

**5. Extension to the case where the distribution of $U$ is not known.** It is possible to generalize the estimator $\hat\beta$ so that it applies to settings where the distribution of $U$ is estimated from data. At least two cases of this type can arise in practice. First, we may observe direct data $U_1,\ldots,U_N$ on $U$, and from those data we may construct an explicit estimator, $\widehat\mu_j(U) = N^{-1}\sum_k (U_i - \bar U)^j$, of $\mu_j(U)$. Here, $\bar U = N^{-1}\sum_k U_k$. Replacing $\nu_j$ by $\widehat\mu_j(U)$ at each appearance in (3), and in all other respects defining $\hat\beta$ as at (4), we obtain a new estimator of $\beta^0 = (\beta_0^0,\ldots,\beta_p^0)^{\mathrm{T}}$. The convergence rate of the new estimator is readily seen to be $O_p\{\min(n,N)^{-1/2}\}$.

Second, and arguably more realistically, we may observe replicated values of $W_j$, so that our dataset is comprised of pairs $(W_{ik}, Y_{ik})$, for $1 \le k \le N_i$ and $1 \le i \le n$, where $W_{ik} = X_i + U_{ik}$ and $Y_{ik} = g(X_i) + \varepsilon_{ik}$. Here, the variables $X_i$, $U_{ik}$ and $\varepsilon_{ik}$ are assumed to be totally independent. In longitudinal data analysis the $N_i$'s are usually small, in the range 2 to 5.

Let us suppose that the distribution of $U$ is symmetric; this would often be a reasonable assumption, and should it fail, a modified version of the argument below could be employed. Let $\mathcal{S}_i$ denote the set of $N_i(N_i - 1)$ distinct pairs $(k_1, k_2)$ with $1 \le k_1 \ne k_2 \le N_i$, and put $N = \sum_{i \le n} N_i(N_i - 1)$. We may estimate the moments $\mu_j(V) = E(V^j)$ of the distribution of $V = U_1 + U_2$, where $U_1$ and $U_2$ are independent copies of $U$,

$$\widehat\mu_j(V) = \frac{1}{N} \sum_{i=1}^n \sum_{(k_1,k_2)\in \mathcal{S}_i} (W_{ik_1} - W_{ik_2})^j.$$

Of course, $\widehat\mu_j(V) = 0$ if $j$ is odd. Let $\mathrm{clt}_{2j}$ and $\mathrm{mnt}_{2j}$ denote the functions that give the $2j$th cumulant, $\kappa_{2j}(Z)$, of a general random variable $Z$ in terms of its moments, and the $2j$th moment in terms of the cumulants,

$$\kappa_{2j}(Z) = \mathrm{clt}_{2j}\{\mu_2(Z),\ldots,\mu_{2j}(Z)\},$$
$$\mu_{2j}(Z) = \mathrm{mnt}_{2j}\{\kappa_2(Z),\ldots,\kappa_{2j}(Z)\}.$$

The $2r$th cumulant of the distribution of $U$ equals half the $2r$th cumulant of the distribution of $V$, and so we define, in succession,

$$\hat\kappa_{2j}(V) = \mathrm{clt}_{2j}\{\widehat\mu_2(V),\ldots,\widehat\mu_{2j}(V)\},$$
$$\widehat\mu_{2j}(U) = \mathrm{mnt}_{2j}\{\tfrac{1}{2}\hat\kappa_2(V),\ldots,\tfrac{1}{2}\hat\kappa_{2j}(V)\},$$

and $\widehat\mu_j(U) = 0$ for odd $j$. Provided the number of indices $i$ in the range $1 \le i \le n$, for which $N_i \ge 2$, increases at rate $n$, the convergence rate of the new estimator is $O_p(n^{-1/2})$.

Next we briefly address hypothesis testing when the distribution of $U$ is not known and it is assumed that (12) holds. We treat in turn the two earlier settings. First, if direct data $U_1,\ldots,U_N$ on $U$ are observed, then we may construct an explicit characteristic-function estimator, $\hat f_U^{\mathrm{Ft}} = n^{-1}\sum_j e^{itU_j}$.

(Here and below, $\hat{f}_U^{\mathrm{Ft}}$ denotes an estimator of $f_U^{\mathrm{Ft}}$, rather than the Fourier transform of an estimator $\hat{f}_U$ of $f_U$.) We replace $f_U^{\mathrm{Ft}}$ in (7) by $|\hat{f}_U^{\mathrm{Ft}}|$, perhaps incorporating a ridge parameter to make the procedure more robust. [Note that, assuming (12), $f_U^{\mathrm{Ft}} = |f_U^{\mathrm{Ft}}|$.] This gives a new version of $\widetilde{F}_X$, leading directly to new formulae for $\bar{F}_X$ and $\hat{F}_X$. In the direct-data setting we do not alter the definitions of $\hat{\psi}(t)$ and $\hat{\phi}(t\,|\,\beta)$, except for replacing $f_U^{\mathrm{Ft}}$ by $\hat{f}_U^{\mathrm{Ft}}$ in the latter.

Second, if (12) holds and we observe replicated data $(W_{jk}, Y_{jk})$, define

$$\hat{f}_U^{\mathrm{Ft}}(t) = \left| \frac{1}{N} \sum_{j=1}^{n} \sum_{(k_1, k_2) \in \mathcal{S}_j} \cos\{t(W_{jk_1} - W_{jk_2})\} \right|^{1/2},$$

potentially incorporating weights to reduce variability. Substituting $\hat{f}_U^{\mathrm{Ft}}$ for $f_U^{\mathrm{Ft}}$ in (7), and modifying $\hat{\psi}(t)$, $\hat{\phi}(t\,|\,\beta)$ and $\hat{\beta}$ by incorporating the replicated data, we obtain an analogue of $T(\hat{\beta})$ which does not require knowledge of $f_U^{\mathrm{Ft}}$. One can also develop analogues, in the case where the distribution of $U$ is estimated from data, of bootstrap methods for calibration.

## APPENDIX: OUTLINE PROOF OF THEOREM 1

Define $p_j = P_j(\nu_0, \ldots, \nu_j)$, $\delta_j = \hat{a}_j - a_j$, $\Delta_j = \hat{A}_j - A_j$,

$$\delta_B^{(k)} = \sum_{\ell=0}^{k} \binom{k}{\ell} \delta_{k-\ell} p_\ell, \qquad \Delta_B^{(k)} = \sum_{\ell=0}^{k} \binom{k}{\ell} \Delta_{k-\ell} p_\ell.$$

Let $\Delta_B$ denote the $(p+1)$-vector with $k$th component $\Delta_B^{(k)}$, and let $\Delta_M$ be the $(p+1) \times (p+1)$ matrix with $(k_1, k_2)$th component $\delta_b^{k_1+k_2}$. Provided (14) holds, the matrix $M$ is finite and strictly positive definite.

In the notation above, $\hat{B} = B + \Delta_B$ and $\widehat{M} = M + \Delta_M$. Therefore, by the Taylor expansion,

(A.1) $\qquad \hat{\beta} = (M + \Delta_M)^{-1}(B + \Delta_B) = \beta^0 + Q + O_p(n^{-1}),$

where $Q = (Q^{(0)}, \ldots, Q^{(p)})^{\mathrm{T}} = M^{-1}\Delta_B - M^{-1}\Delta_M M^{-1} B$ is a $(p+1)$-vector. Since $Q$ is expressible exactly as the mean of $n$ independent and identically distributed random $(p+1)$-vectors with zero expected value, it is readily proved that $n^{1/2}Q$ is asymptotically normally distributed with zero mean and finite variance.

With $W_j = X_j + U_j$ and $Y_j = g(X_j) + \varepsilon_j$ denoting the data, we have

(A.2) $\quad Q^{(k)} = \sum_{\ell=0}^{p} (M^{-1})_{k\ell} \Delta_B^{(\ell)} - \sum_{\ell_1=0}^{p} \sum_{\ell_2=0}^{p} (M^{-1})_{k\ell_1} (\Delta_M)_{\ell_1 \ell_2} (M^{-1}B)_{\ell_2}$

$$= \sum_{\ell=0}^{p}\sum_{r=0}^{\ell}\binom{\ell}{r}(M^{-1})_{k\ell}p_r\frac{1}{n}\sum_{j=1}^{n}(1-E)\left\{\sum_{s=0}^{p}\beta_s^0 X_j^s W_j^{l-r}+\varepsilon_j W_j^{l-r}\right\}$$

$$(A.3) \qquad -\sum_{\ell_1=0}^{p}\sum_{\ell_2=0}^{p}\sum_{r=0}^{\ell_1+\ell_2}\binom{\ell_1+\ell_2}{r}(M^{-1})_{k\ell_1}(M^{-1}B)_{\ell_2}p_r$$

$$\times \frac{1}{n}\sum_{j=1}^{n}(1-E)W_j^{l_1+l_2-r},$$

where $E$ denotes the expectation operator. Note too that

$$\hat{\psi}(t)f_U^{\text{Ft}}(t) = \frac{1}{n}\sum_{j=1}^{n}Y_j e^{itW_j} = \sum_{k=0}^{p}\beta_k^0\frac{1}{n}\sum_{j=1}^{n}X_j^k e^{itW_j} + \frac{1}{n}\sum_{j=1}^{n}\varepsilon_j e^{itW_j},$$

$$\hat{\phi}(t\,|\,\beta)f_U^{\text{Ft}}(t) = \sum_{k=0}^{p}\beta_k\sum_{\ell=0}^{k}\binom{k}{\ell}\left(\frac{1}{n}\sum_{j=1}^{n}W_j^\ell e^{itW_j}\right)\phi_{k-\ell}(t),$$

where $\phi_r(t) = f_U^{\text{Ft}}(t)(i^{-1}D_t)^r f_U^{\text{Ft}}(t)^{-1}$. Define

$$\chi_{1k}(t) = \frac{1}{n}\sum_{j=1}^{n}\left\{X_j^k - \sum_{\ell=0}^{k}\binom{k}{\ell}W_j^\ell \phi_{k-\ell}(t)\right\}e^{itW_j},$$

$$(A.4) \qquad \chi_2(t) = \frac{1}{n}\sum_{j=1}^{n}\varepsilon_j e^{itW_j},$$

$$\chi_{3k}(t) = \sum_{\ell=0}^{k}\binom{k}{\ell}\left(\frac{1}{n}\sum_{j=0}^{n}W_j^\ell e^{itW_j}\right)\phi_{k-\ell}(t).$$

In this notation,

$$\{\hat{\psi}(t)-\hat{\phi}(t\,|\,\hat{\beta})\}f_U^{\text{Ft}}(t) = \sum_{k=0}^{p}\beta_k^0\chi_{1k}(t) + \chi_2(t)$$

$$-\sum_{k=0}^{p}(\hat{\beta}_k-\beta_k^0)\sum_{\ell=0}^{k}\binom{k}{\ell}\left(\frac{1}{n}\sum_{j=1}^{n}W_j^\ell e^{itW_j}\right)\phi_{k-\ell}(t).$$

The series multiplying $(\hat{\beta}_k-\beta_k^0)$ in the last term equals $\chi_{3k}(t)$. Using (A.1),

$$\left\{\int|\hat{\psi}(t)-\hat{\phi}(t\,|\,\hat{\beta})|^2 f_U^{\text{Ft}}(t)^2 w(t)\,dt\right\}^{1/2}$$

$$(A.5) \qquad = \left\{\int\left|\sum_{k=0}^{p}\beta_k^0\chi_{1k}(t)+\chi_2(t)-\sum_{k=0}^{p}(\hat{\beta}_k-\beta_k^0)\chi_{3k}(t)\right|^2 w(t)\,dt\right\}^{1/2}$$

$$+ O_p(n^{-1})$$
$$= \xi_n^{1/2} + O_p(n^{-1}),$$

where

$$\xi_n = \int \left| \sum_{k=0}^{p} \beta_k^0 \chi_{1k}(t) + \chi_2(t) - \sum_{k=0}^{p} Q^{(k)} \chi_{3k}(t) \right|^2 w(t)\, dt$$

and $Q^{(k)}$ is as in (A.3).

Note that $\chi_{3k}(t) \to \xi_{3k}(t)$ as $n \to \infty$, where

$$\xi_{3k}(t) = \sum_{\ell=0}^{k} \binom{k}{\ell} E(W^\ell e^{itW}) \phi_{k-\ell}(t) = f_U^{\mathrm{Ft}}(t)(i^{-1} D_t)^k f_X^{\mathrm{Ft}}(t).$$

Using standard properties of sums of independent random variables, and referring to (A.3)–(A.4) to deduce the relationships among $\chi_{1k}(t)$, $\chi_2(t)$ and $Q^{(\ell)}$, it may be proved that

(A.6) $$n\xi_n \to \xi$$

in distribution as $n \to \infty$, where

$$\xi = \int \left| \sum_{k=0}^{p} \beta_k^0 \xi_{1k}(t) + \xi_2(t) - \sum_{k=0}^{p} R^{(k)} \xi_{3k}(t) \right|^2 w(t)\, dt,$$

$\xi_{10}, \ldots, \xi_{1p}$, $\xi_2$ and $R^{(1)}, \ldots, R^{(p)}$ are jointly distributed, $\xi_{1k}(t)$ and $\xi_2(t)$ are complex-valued Gaussian processes with zero means, $R = (R^{(1)}, \ldots, R^{(p)})^{\mathrm{T}}$ is a Gaussian $(p+1)$-vector with zero mean, and the covariances among $\xi_{1k}(t)$, $\xi_2(t)$ and $R^{(u)}$ are identical to the covariances among

$$\left\{ X^k - \sum_{\ell=0}^{k} \binom{k}{\ell} W^\ell \phi_{k-\ell}(t) \right\} e^{itW}, \qquad \varepsilon e^{itW}$$

and

$$\sum_{\ell=0}^{p} \sum_{r=0}^{\ell} \binom{\ell}{r} (M^{-1})_{ul} p_r \left\{ \sum_{s=0}^{p} \beta_s^0 (1-E) X^s W^{l-r} + \varepsilon W^{l-r} \right\}$$

$$- \sum_{\ell_1=0}^{p} \sum_{\ell_2=0}^{p} \sum_{r=0}^{\ell_1+\ell_2} \binom{\ell_1+\ell_2}{r} (M^{-1})_{ul_1} (M^{-1} B)_{\ell_2} p_r (1-E) W^{l_1+l_2-r},$$

respectively. Note particularly that these covariances depend on the distribution of $\varepsilon$ only through the variance of that quantity.

Part (a) of Theorem 1 follows from (A.5) and (A.6). The fact that $P(\xi > 0) = 1$ can be deduced from the fact that $\mathrm{var}\,\varepsilon > 0$. Derivation of part (b) is virtually identical to that of part (a). To prove part (c) of Theorem 1, note that the presence of the perturbation $n^{-1/2}c\gamma(x)$ in (17) influences $Q^{(k)}$, at (A.2), only by adding $n^{-1/2}c\eta_{1k}$ to the first term on the right-hand side of (A.2) [and, hence, by adding the same quantity to the far right-hand side at (A.3)], where

$$\eta_{1k} = \sum_{\ell=0}^{p}\sum_{r=0}^{\ell} \binom{\ell}{r} (M^{-1})_{k\ell} E\{\gamma(X)W^{\ell-r}\}p_r.$$

The impact of the perturbation on $\hat{\psi}(t)$ can be described completely by adding $n^{-1/2}c\eta_2(t)$ to $\chi_2(t)$, where $\eta_2(t) = E\{\gamma(X)e^{itX}\}f_U^{\text{Ft}}(t) = E\{\gamma(X)e^{itW}\}$. The perturbation has no effect on $\hat{\phi}(\cdot\,|\,\beta)$.

Therefore, retracing the arguments leading to (A.5), we see that result continues to hold if we add, within the modulus signs in the definition of $\xi_n$, the quantity $n^{-1/2}c\eta_3(t)$, where $\eta_3(t) = \eta_2(t) - \eta_4(t)$ and $\eta_4(t) = \sum_{0\leq k\leq p}\eta_{1k}\chi_{3k}(t)$. It follows from this property that part (c) of Theorem 1 holds provided $\eta_3$ is nonzero on a set of positive measure. In the next paragraph we derive this property.

The constants $\eta_{10},\ldots,\eta_{1p}$ are the unique solutions of the equations

$$E\{\gamma(X)W^k\} = E\left\{\left(\sum_{j=0}^{p}\eta_{1j}X^j\right)W^k\right\}, \qquad 0\leq k\leq p.$$

Equivalently, $\gamma_1(x) = \sum_{0\leq j\leq p}\eta_{1j}x^j$ is the unique $p$th degree polynomial for which $E\{\gamma(X)W^k\} = E\{\gamma_1(X)W^k\}$ for $0\leq k\leq p$. From this property, and the fact that $f_U^{\text{Ft}}$ vanishes only at isolated points and

$$\eta_3(t) = E\{\gamma(X)e^{itW}\} - E\{\gamma_1(X)e^{itW}\}$$
$$= E[\{\gamma(X) - \gamma_1(X)\}e^{itX}]f_U^{\text{Ft}}(t),$$

we deduce that $\eta_3$ vanishes almost everywhere if and only if $\gamma = \gamma_1$ almost everywhere on the support of the distribution of $X$. However, the conditions imposed for part (c) of Theorem 1 rule this out, and so $\eta_3$ is nonzero on a set of positive measure.

20 P. HALL AND Y. MA[27] STEFANSKI, L. A. and CARROLL, R. J. (1987). Conditional scores and optimal scores for generalized linear measurement-error models. *Biometrika* **74** 703–716. MR0919838
[28] STEFANSKI, L. A. and CARROLL, R. J. (1990). Deconvoluting kernel density estimators. *Statistics* **21** 169–184. MR1054861
[29] TAUPIN, M. (2001). Semi-parametric estimation in the nonlinear structural errors-in-variables model. *Ann. Statist.* **29** 66–93. MR1833959
[30] TSIATIS, A. A. and MA, Y. (2004). Locally efficient semiparametric estimators for functional measurement error models. *Biometrika* **91** 835–848. MR2126036

DEPARTMENT OF MATHEMATICS AND STATISTICS
UNIVERSITY OF MELBOURNE
PARKVILLE, VICTORIA 3010
AUSTRALIA
E-MAIL: hall@unimelb.edu.au

DEPARTMENT OF STATISTICS
TEXAS A&M UNIVERSITY
3143 TAMU
COLLEGE STATION, TEXAS 77843-3143
USA
E-MAIL: ma@stat.tamu.edu